\documentclass[11pt]{article}
\usepackage{amssymb,amsmath,enumerate,latexsym}

\usepackage{float}

\usepackage{graphicx}

\def\mobius{M\"obius }
\newtheorem{lemma}{Lemma}
\newtheorem{proposition}{Proposition}
\newtheorem{theorem}{Theorem}

\newtheorem{corollary}{Corollary}

\newtheorem{assertion}{Assertion}
\newtheorem{question}{Question}

\newenvironment{proof}{\trivlist
\item[\hskip\labelsep{\it Proof}\,:]}{\hfill{$\Box$}\endtrivlist}

\newcommand{\re}{\mathop{\rm Re}\nolimits}
\newcommand{\im}{\mathop{\rm Im}\nolimits}

\def\R{\mathbb{R}}
\def\N{\mathbb{N}}
\def\C{\mathbb{C}}
\def\Z{\mathbb{Z}}
\def\esf{\mathbb{S}}

\def\lc{\mathcal{L}}

\def\ve{{\varepsilon}}

\def\rth{{\mathbb{R}^3}}

\begin{title}
{Bending the Helicoid}
\end{title}

\begin{author}
{William H. Meeks III\thanks{This material is based upon
 work for the NSF under Award No. DMS - 0405836.} 
 \and  Matthias Weber\thanks{This material is based upon
 work for the NSF under Award No. DMS - 0139476 and DMS - 0505557. Any opinions, findings, and conclusions or recommendations
 expressed in this publication are those of the authors and do not
 necessarily reflect the views of the NSF}}
\end{author}
\begin{document}

\maketitle
\begin{abstract}

  We construct Colding-Minicozzi limit minimal laminations in open
 domains in $\rth$ with the singular set of $C^1$-convergence
  being any properly embedded $C^{1,1}$-curve. By
 Meeks' $C^{1,1}$-regularity theorem, the singular set of
   convergence of a Colding-Minicozzi limit minimal lamination
 ${\cal L}$ is a locally finite collection
    $S({\cal L})$ of $C^{1,1}$-curves that are orthogonal to 
the leaves of the lamination. Thus, our existence theorem gives a complete answer
 as to which curves appear as
     the singular set of a Colding-Minicozzi limit minimal lamination.
 
  In the case the curve is the unit circle $\esf^1(1)$ in the $(x_1, x_2)$-plane,
      the classical Bj\"orling  theorem produces an infinite sequence of complete minimal
 annuli $H_n$ of finite total curvature
     which contain the circle. The complete minimal surfaces $H_n$ 
contain embedded compact
 minimal annuli $\overline{H}_n$ in closed compact neighborhoods $N_n$ of 
the circle that converge as $n \to \infty$ to $\rth - x_3$-axis. In this case, we prove that the $\overline{H}_n$ converge  on compact sets to the foliation of $\rth - x_3$-axis by vertical half planes with boundary the $x_3$-axis and with $\esf^1(1)$ as the singular set of $C^1$-convergence.  The $\overline{H}_n$
      have the appearance of highly spinning helicoids with the circle
 as their axis and are named {\em bent helicoids}.

\vspace{.17cm}


\end{abstract}

\section{Introduction.}

In \cite{cm23}, Colding and Minicozzi consider the question of the 
compactness of a sequence $\{M_{n}\}_{n\in\N}$ of embedded minimal surfaces 
in a
Riemannian three-manifold $N$ which are {\it locally simply
connected} in the following sense: for each small open geodesic ball in $N$
 and for each $n$ sufficiently large,
 $M_{n}$ intersects the ball in disk components, with each disk component having its boundary in the boundary of the ball.  They prove
 that every
such sequence of minimal surfaces has a subsequence which converges
to a possibly singular limit minimal lamination $\lc$ of $N$. In
certain cases, the minimal lamination $\lc$ is nonsingular and is a
minimal foliation of $N$. In this case, they prove that the singular set of
$C^1$-convergence consists of a properly embedded locally finite
collection $S(\lc)$ of Lipschitz curves that intersect the leaves of $\lc$
transversely; we call such a limit foliation $\lc$ a {\it
Colding-Minicozzi limit minimal lamination}.

In \cite{mr8}, Meeks and Rosenberg applied these results of
Colding-Minicozzi to prove that the plane and the helicoid are the
only properly embedded simply connected minimal surfaces in $\rth$.
A standard blow-up argument then shows that small neighborhoods of points
of large almost-maximal curvature on an embedded minimal surface of positive 
injectivity radius in a homogeneously regular three-manifold $N$ have the 
appearance of homothetically shrunk helicoids. 
 An  application by Meeks \cite{me30,
me25} of this local picture for a minimal disk centered at a point
of large almost-maximal curvature demonstrates
  that the singular curves
$S(\lc)$ of a Colding-Minicozzi lamination $\lc$ have class $C^{1, 1}$ and are
orthogonal to the leaves of $\lc$.  The proof by Meeks of
the $C^{1,1}$-regularity of $S(\lc)$ leads naturally to a unique
related lamination metric on the minimal foliation 
on $\lc$ of $N$ (\cite{me30}).  This regularity
theorem and lamination metric theorem have useful applications which include
the classification of properly embedded minimal surfaces of finite
genus in $\rth$ and in other three-manifolds (see \cite{mpr3, mpr4,
mpr8, mpr10, mr13, mr10}).

In all previously considered examples of sequences of locally simply
connected minimal surfaces which converge to a
minimal foliation $\lc$ with nonempty singular set of $C^1$-convergence $S(\lc)$, 
the set $S(\lc)$
consisted of geodesics.  While the first author had thought that this property
might hold in general, it was pointed out to him by Frank Morgan that
 it was reasonable to expect that there exist compact embedded minimal annuli
$A_n$ that would converge to a Colding-Minicozzi limit minimal lamination
 $\widetilde{\lc}$ of an open set of $\rth$ and with $S(\widetilde{\lc})$ being the unit circle in the 
$(x_1, x_2)$-plane. In fact, the following main theorem shows that  in the case of $\rth$, {\em any} 
 $C^{1,1}$-curve $S(\lc)$ occurs as a singular set of $C^{1}$-convergence of a 
Colding-Minicozzi limit minimal lamination.

\begin{theorem} \label{thm1}   Every properly embedded $C^{1,1}$-curve
   $\alpha$ in an open set $O$ in $\rth$ has a neighborhood
   foliated by a particular Colding-Minicozzi limit minimal lamination $\lc$
 with singular set
   of $C^1$-convergence being $\alpha$. The minimal leaves of
 this lamination $\lc$
 are a $C^{1,1}$-family of pairwise disjoint flat disks of varying radii.
The disks are centered
 along and orthogonal to $\alpha$. More generally, if $N$ is a
 closed regular neighborhood of $\alpha$ formed by disjoint flat disks
 orthogonal to $\alpha$ and $N'$ is a similarly defined foliation in the
 interior of $N$, then $N'$ is contained in a Colding-Minicozzi minimal
 lamination which lies in $N$.
\end{theorem}

The main step in the proof of Theorem~\ref{thm1} is to first prove the
 theorem when $\alpha$ is analytic with a compact exhaustion
 $\alpha(1) \subset \alpha(2) \subset ... \subset \alpha(n) \subset ...$,
 where $\alpha (i)$ is a compact connected arc in $\alpha$. We do this by 
giving an essentially explicit construction of a
sequence of embedded compact {\em bent helicoids} $H_{\alpha,n}$ which contain
$\alpha(n) \subset \alpha$ as an ``axis'' and whose Gauss maps rotate
faster and faster along $\alpha(n)$ as $n
\rightarrow \infty$.  In this case, the $H_{\alpha,n}$ converge to a family 
of pairwise disjoint flat disks of varying radii orthogonal to $\alpha$. 
The construction of the $H_{\alpha,n}$ is based on the
classical Bj\"orling formula.  Our main difficulty in proving Theorem
\ref{thm1} in the analytic case is to demonstrate the embeddedness of the
 $H_{\alpha,n}$ in a {\em fixed}
neighborhood of $\alpha(n)$. The general case of the theorem follows from
 the analytic case by approximating $\alpha$ by a sequence of embedded
 analytic curves with uniformly locally bounded curvature, which is always
 possible for $C^{1,1}$-curves.  

In the special case that $\alpha$ is the
unit circle in the $(x_1, x_2)$-plane, then, for all $n \in \N$, we can choose
 $\alpha(n) = \alpha$ and each compact annular bent helicoid 
$\overline{H}_n = H_{\alpha,n}$ contains
$\alpha$ and is the image of a compact portion of a globally defined 
explicit periodic 
complete minimal immersion 
$f_n \colon \mathbb{C} \rightarrow \rth$. In this case, we let $H_n$
 denote the image complete minimal annulus
 $f_n(\C)$ and define compact
 embedded annuli $\overline{H}_n \subset H_n$ which converge 
to the limit minimal
 foliation $\lc$ of $\rth-x_3$-axis by vertical half
planes and with $S(\lc) = \alpha$. We refer the reader to Section 3
for the analytic description of the parametrizations $f_n$ of these
 special bent helicoids whose coordinate functions are expressed 
in terms of real and imaginary parts of the 
 functions $\cos(z)$ and $\sin(z)$ for $z \in \C$. We 
also describe the analytic
 Weierstrass data for their image finite total curvature annuli $H_n$ 
in terms of simple rational functions on the punctured complex plane 
$\mathbb{C} - \{0 \}$.

The complete minimal annulus $H_n$ has finite total curvature
 $-4 \pi(n+1)$ with the dihedral group $D(2n)$ of symmetries and 
contains $n$ lines in the $(x_1, x_2)$-plane passing through the
 origin. The large symmetry group and the explicit representation
 of $H_n$ allows us to define the compact {\it embedded} annuli
 $\overline{H}_n \subset H_n$ which converge to the minimal foliation
 $\lc$ of $\rth$. By way of approximation, this special
 case of a circle plays a key role in the proof of Theorem~\ref{thm1} 
in the more general
 case where $\alpha$ is an arbitrary properly embedded analytic
 curve in an open set $O$. This is because at every point of the
 analytic curve the related bent helicoids that we construct are 
closely approximated by the related bent helicoids of the second order
 approximately osculating circle at the point. 
Based on our construction of these bent helicoids and the $C^{1,1}$-regularity 
theorem of Meeks \cite{me30, me25}, we ask the following
related question.

\begin{question} \label{quest1} Is there a natural generalization of
 Theorem~\ref{thm1} to Riemannian three-manifolds?
\end{question}

\begin{figure}[H]
   \centerline{\includegraphics[width=4in]{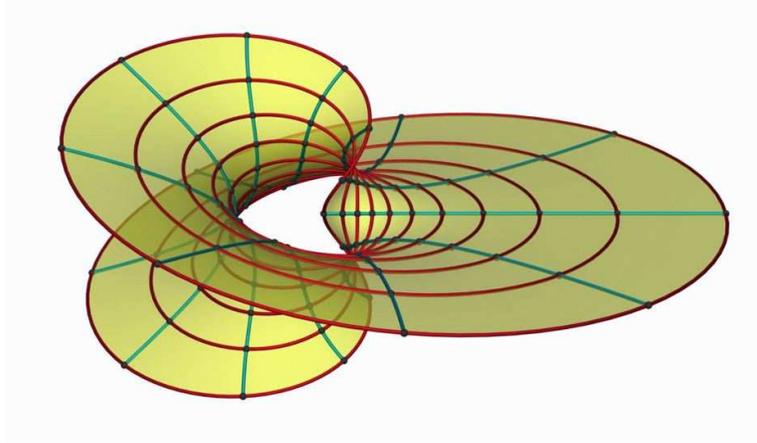}}
   \caption{The complete minimal \mobius strip}
   \label{fig:mobius}
\end{figure}

It turns out that the bent helicoids $H_n$ also make sense for values $n= k-\frac{1}{2}$, where $k \in \N$, and for these fractional values the image surface is a complete immersed minimal \mobius strip. One special case of these bent helicoids was known before, namely
 $n=\frac12$. This example is the
Meeks' Minimal \mobius strip of finite total curvature $-6 \pi$ defined in \cite{me7}, see
 Figure \ref{fig:mobius}.

For larger integer values $n$, near the unit circle the surface
 $H_n$ looks  like a bent
helicoid:

\begin{figure}[H]
   \centerline{\includegraphics[width=3.8in]{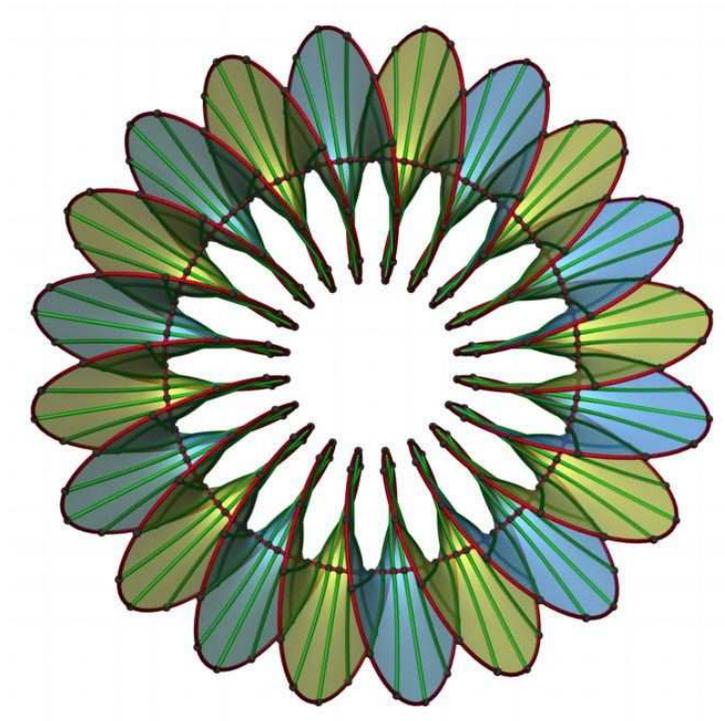}}
   \caption{The bent helicoid for $n$=10}
   \label{fig:bentheli10}
\end{figure}

We would like to thank Bruce Solomon for helpful conversations about regularity questions.

\section{Bj\"orling's theorem and the analytic representation of bent
 helicoids in the circular case.}

We now recall Bj\"orling's theorem \cite{dhkw}.
Let $c \colon [a,b]\to \R^3$ be any real analytic
 curve  and
$n \colon [a,b]\to \esf^2 \subset \R^3$ be  any real analytic vector field 
 perpendicular to $c'(t)$.
 Consider $[a,b]\times\{0\} \subset \C$.
 By analyticity, there are
  a small positive $\ve$ and   unique 
holomorphic extensions
$c \colon [a,b]\times (-\ve, \ve)\to \C^3$, and 
$n:[a,b]\times (-\varepsilon, \varepsilon)\to \C^3$.
Using these extensions, we define for $z = x+yi \in \C$
\[
 F (z) = \re\left( c(z) -  i \int_0^z n(w)\times c'(w)\, dw \right).
\]
This is a minimal map that extends $c$  and $n$ in the sense that 
for  $t\in[a,b]$, $F(t)=c(t)$ and $n(t)$ is the surface normal.

This formula was used by H.A.~Schwarz to prove the classical reflection principles 
for minimal surfaces and to prove that the
helicoid is the only ruled minimal surface besides the plane. 

On the other hand, the above formula has not produced other globally
 interesting examples of minimal surfaces. This is mainly due
 to the fact that the Bj\"orling 
integral is usually impossible to evaluate explicitly, making
 it hard to say something about global properties of the minimal surfaces.
In the case of the unit circle in the $(x_1, x_2)$-plane, there is a 
natural analytic parametrization as well as a natural sequence of analytic unit normal fields for which we can explicitly evaluate this integral to produce a sequence of Bj\"orling surfaces.

Let $c(t)=(\cos t, \sin t,0)$ be the usual parametrization of the unit circle in the $(x_1,x_2)$-plane. A
basis of normal fields along $c$ is given by
\begin{align*}
n_1(t)&=-c(t),
\\
n_2(t)&=(0,0,1).
\end{align*}

Define a new normal field with relative rotational speed $a \in \R^+$ by

$$
n(t)=\cos(a t)n_1(t) + \sin(at) n_2(t).
$$
For $z \in \mathbb{C}$, we let $c(z), n(z)$ be the related vector
 valued holomorphic functions mapping $\mathbb{C}$ to $\mathbb{C}^3$.
Then, using Bj\"orling's formula \cite{dhkw},

\[
F (z) = \re\left( c(z) -  i \int_0^z n(w)\times c'(w)\, dw
\right)
\]
defines a minimal surface with parameter domain $\C$ which extends the 
circle $c(t)$ and has $n(t)$ as the
Gauss map along the circle. We refer to this surface as the {\em bent
helicoid} $H_a$. Here we consider $c \colon \R \to \rth$ to be a parametrized curve with related $n(t)$ along it; when $a \in \N \subset \R$, then $n(t)$ is well defined on the image circle $\esf^{1}(1)$.

For a fixed value $a \in \R^+$, the conformal harmonic map $F \colon \C \to \rth$ is
 explicitly calculated to be:
\begin{align*}
F(z)&= \re\int_0^z\begin{pmatrix} i \cos (w) \sin (a w)-\sin (w)
\\
\cos (w)+i \sin (w) \sin
    (a w)
  \\
  i \cos (a w)
  \end{pmatrix}\,dw
\\
& \overset{a\ne1}{=}
\re
\begin{pmatrix}
\cos (z)-\frac{i (\cos (z) \cos (a z) a-a+\sin (z)
    \sin (a z))}{a^2-1}
    \\
    \sin (z)-\frac{i (a \cos (a z) \sin
    (z)-\cos (z) \sin (a z))}{a^2-1}
    \\
\frac{i \sin (a
    z)}{a}
    \end{pmatrix}.
\end{align*}

Moreover, one can now convert this data for $H_a$ to data for the
 classical Weierstrass representation
\[
 F(z) = \frac12 \int^z \left(  \frac1G -G, i \left(\frac1G +G\right) ,
 2\right)\cdot dh.
\]
This conversion produces a stereographically projected Gauss map
$$G(z)  = -e^{i z} \frac{\cos(az)}{1-\sin(az)}$$
and a complexified height differential
$$
dh = i \cos(az) dz.
$$

After the substitution $w=e^{i z},$ and for $a = n \in \N$, then 
\begin{align*}
G(w)&=-w \frac{ w^n+i}{i w^n+1},
\\
dh&=\frac{1}{2w}(w^n+w^{-n})\, dw.
\end{align*}
Thus, in this case, we see that $H_a = H_n$ is a complete minimal surface of a finite total
curvature $-4\pi (n+1)$ and with parameter domain $\C - \{ 0 \}$.

For $n=0$, we recover the familiar Weierstrass representation of the catenoid.

\section{The geometry and embeddedness of fundamental pieces of bent helicoids in the circular case}

We now collect some simple properties of the bent helicoids $H_a$:

\begin{proposition}
\begin{enumerate}
\item For $a \in \N$, the immersion
 $F(z) \colon \mathbb{C} \to \rth$ is $2\pi$ periodic in the sense that $F(z) = F(z + 2\pi)$.
\item
The vertical coordinate  lines $x =t_k=\frac{2k+1}{2a}\pi$ 
are mapped to the straight lines
$s \mapsto s (\cos(t_k), \sin(t_k),0)$
\item
The $180$ degree rotations around the points $t_k=\frac k{a} \pi$ 
in the domain $\C$ induce isometries of the surface which are 180 degree 
rotations about the lines $s \mapsto s (\cos(t_k), \sin(t_k),0)$ (orthogonal to the surface)
in $\rth$.
\item
The surface is invariant under rotation by angle $\frac{\pi}{a}$ about
 the $x_3$-axis.
\end{enumerate}
\end{proposition}

\begin{proof}
The first claim is trivial. We compute $f(t_k+t i)$ to be

$$
\frac{(-1)^k \cosh (a t) \sinh (t)+\cosh (t) \left(a^2-(-1)^k \sinh (a t)
   a-1\right)}{a^2-1}
   \left(
\cos \left(t_k\right),\sin \left(t_k\right),0
   \right),
   $$
  which proves the second claim. Alternatively, one can also see this from the  uniqueness of the Bj\"orling solution as
 follows.
At the points $t_k = \frac{2k+1}{2a}\pi$, we have
\begin{align*}
c(t_k)&=(\cos(t_k), \sin(t_k),0),
\\
n(t_k) &= (0,0,(-1)^k).\end{align*}
 Since a rotation around the line $s \mapsto s (\cos(t_k), \sin(t_k),0)$
 maps the Bj\"orling data $c$ and $n$ to $c$ and $-n$,
 it must map the surface $H(a)$ to the same surface with reversed orientation
 (by the uniqueness of the Bj\"orling solution). Because the line is tangent
 to the surface at $c(t_k)$, it must lie entirely on the surface.

The remaining claims are proven in a similar fashion.
\end{proof}

\begin{proposition}
\label{prop:conformal}
The conformal factor of the metric of $H(a)$ with conformal parametrization
 $F(x,y) = F(x+yi) \colon
\mathbb{C} \to \rth$ is
$$
\lambda(x,y)=|F_x|=|F_y|=\cosh (y) \cosh (a y)-\sin (a x) \sinh (y)
$$
and the tangent vector of the curve $y\mapsto F(0,y)$ is
$$
F_y(0,y)=(\sinh (y),\sinh (y) \sinh (a y),-\cosh (a y)).
$$
\end{proposition}
\begin{proof}By direct computation.
\end{proof}

Moreover, the simple form of the Weierstrass data allows us to
 establish some other remarkable properties
of this family of surfaces, which we now describe.

The following lemma shows that the image of the vertical half-lines $T\mapsto x\pm T i$ are (for large $T$) close to horizontal half-lines in space.

\begin{lemma}
\label{lem:arc}
\begin{align*}
\lim_{T\to\infty} e^{-(a+1)T} F(x+T i) &= \frac{1}{4(a+1)} \left( -\sin((a+1)x),
\cos((a+1)x),0 \right),
\\
\lim_{T\to -\infty} e^{(a+1)T} F(x+T i) &= \frac{1}{4(a-1)} \left( \sin((a+1)x),
-\cos((a+1)x),0 \right).
\end{align*}
\end{lemma}
\begin{proof}
This follow from the integrated formula of $F$ by straightforward computation.
\end{proof}

\begin{corollary}
For $a>2$, the image under $F$ of $[-\frac{\pi}{2a},\frac{\pi}{2a}]\times 
(-\infty,\infty)$ is embedded, except at the origin in $\rth$ where the two 
boundary lines intersect. Hence, the image under $F$ of the fundamental piece $(\frac{-\pi}{2a}, \frac{\pi}{2a} ] \times
(-\infty, \infty)$ is an embedded surface. 
\end{corollary}
\begin{proof}
Subdivide for any $T>0$  the domain $[-\frac{\pi}{2a}, \frac{\pi}{2a} ] \times (-T, T)$ into two closed pieces $R_T^{\pm}$ depending on the sign of
  $y$. The boundary of $F(R_T^{+})$ consists of four pieces: the circular arc $\alpha_1 = \alpha([-\frac{\pi}{2a}, \frac{\pi}{2a}])$, the rays $\alpha_2 = F(\{- \frac{\pi}{2a} \} \times [0,T))$ and 
  $\alpha_3 = F(\{+ \frac{\pi}{2a} \} \times [0, T))$
  and the image arc
$\alpha_4 = F([-\frac{\pi}{2a}, \frac{\pi}{2a}] \times \{T\}$. 

For large $T$, the total curvature of the boundary of $F(R_T^{+})$ is asymptotically $\frac{2\pi}a+3\pi<4\pi$ when $a>2$. In addition, this boundary is embedded, because the respective arcs are disjoint and embedded as individuals.
Thus the entire image $F(R_T^{+})$ is embedded (for $a>2$) by \cite{eww1}  for any large $T$, and thus also $F(R_\infty^{+})$.

As the third coordinate function $\frac1a \cos(ax)\sinh(ay)$ changes
 sign with $y$, $F(R_\infty^{+})$ is contained in the closed upper half space and intersects the $x_1x_2$-plane only in the two rays $\alpha_2$ and $\alpha_3$. The $180^\circ$ rotation about the $x_1$-axis rotates $F(R_\infty^{+})$ into $F(R_\infty^{-})$. Hence the whole surface is embedded except at the origin as the intersection point of the two lines. This proves the corollary.
\end{proof}

\begin{figure}[H]
  \centerline{\includegraphics[width=2.8in]{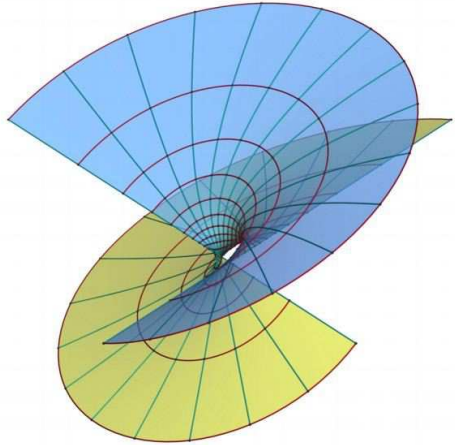}}
  \caption{An embedded fundamental piece of $H_a$ for $a=2$.}
  \label{fig:embedded}
\end{figure}

\section{Approximation results}

This section has two goals. First, in the circle case, we give an explicit estimate for
 how close the $y$-curves $y\mapsto F(x,y)$ are to the lines $L_x$ in $\rth$ passing through $F(x,0)$ and tangent to the curve at this point. This estimate shows that on compact subsets of $\rth$, as $a \to \infty$, the $y$-curves converge $C^k$ to the corresponding line in $ \{ L_x \}_{x \in \R}$, and so, on a given large compact set of $\rth$ and for $a \in \N$ sufficiently large, $H_a$ closely approximates the related ruled surface along the circle. However, our estimates are not sufficient to prove Theorem~\ref{thm1} by comparison with the ruled surface; we get around this problem in the next section by using the large symmetry group of $H_n$. Second, for an analytic curve $\alpha$, we  compare a general bent helicoid $H_{\alpha,n}$ (to be defined)
to the circular helicoid $\overline{H}_a$ in terms of how far the surfaces
 are apart along the related $y$-curve lines that begin near a point $\alpha(x)$, 
when the circle is the second order approximation to $\alpha$ at
 $\alpha(x)$.

For the first part, we compare the minimal surface $H_a$ with a suitably
 parameterized ruled surface
\[
R(x,y)=c(x)+ t_a(y) c'(x)\times n(x),
\]
where
\begin{align*}
t_a(y)&=
\frac{a \cosh (y) \sinh (a y)-\cosh (a y) 
\sinh (y)}{a^2-1}-\sin (a x)
(\cosh (y)-1)
\\
&\approx
\frac{\cosh(y)\sinh(a y)}{a} \qquad\text{for $a$ large}
\end{align*}
is the primitive of $|F_y(x,y)|$.

This ruled surface has the same core circle as $F$ and its
 ruling lines point in the same direction as the tangent vectors
 $F_y(x,0)$. The ruling lines are parameterized so that their
 speed is equal to the speed $|F_y(x,y)|$ given by
 Proposition  \ref{prop:conformal}.
While the ruled surface $R$ is a poor approximation for $F$ when
 $y=\im z$ is large (by Lemma \ref{lem:arc}),
it becomes better and better  in the  range $|y|<d=d(a)$, where $d$ is determined so that the curves $F(x\pm d i)$ stays at distance $1$ away from the core circle
by using the approximate expression for $t_a(y)$ to solve $t_a(d)=1$ for $d$. This motivates (for $a$ large) our definition of $d$:
\[
d=d(a)=\frac{\log(a)}a.
\]
Observe that
\[
|R(x, d)-R(x,0)| \approx 1 \approx |R(x, -d)-R(x,0)|.
\]
This means that the curves $\im z=\pm d$ are approximately
 mapped onto the boundary of the tube of radius 1 around the  unit circle.

Observe also that for $x= t_k=\frac k{a} \pi$, the
 parameterizations $F(t_k,y)$ and $R(t_k,y)$ coincide.
 
 Now we can state and prove our approximation theorem.

\begin{lemma}
For $|y|\le d(a)$, we have
\[
|R(x,y)-F(x,y)|\le d(a)  \qquad\text{and $d(a)\to 0$ as $a\to\infty$}.
\]
\end{lemma}
\begin{proof}
By the definitions, we have
\[
\frac{\partial}{\partial y} \left(R(x,y)-F(x,y) \right)
=
\begin{pmatrix}
\cos (a x) \sinh (y) (\sin (x) \sinh (a y)-\cos (x) \cos (a x))
\\
-\cos (a x) \sinh
   (y) (\cos (a x) \sin (x)+\cos (x) \sinh (a y))
   \\
-\cos (a x) ((\cosh (y)-1) \cosh (a
   y)-\sin (a x) \sinh (y))
\end{pmatrix}.
\]

\noindent By direct computation, we obtain
\[
\left|
\frac{\partial}{\partial y} \left( R(x,y)-F(x,y)\right)
\right|^2
=
4 \cos ^2(a x) \cosh (a y) \sinh ^2\left(\frac{y}{2}\right)
 (\cosh (y) \cosh (a
   y)-\sin (a x) \sinh (y)).
\]
From these formulas, it follows that for $|y|<d$, then
\[
\left|
\frac{\partial}{\partial y} \left( R(x,y)-F(x,y)\right)
\right|
\le
\cosh (a y) |\sinh (y)|,
\]
and so,
\begin{align*}
|R(x,y)-F(x,y)|&\le \int_0^d \cosh (a y) |\sinh (y)|\, dy
\\
&\le
\frac{-\cosh (d) \cosh (a d)+a \sinh (d) \sinh (a d)+1}{a^2-1}
\\
&\le d,
\end{align*}
as claimed.
\end{proof}

Our next goal is to compare the Bj\"orling surfaces associated to arbitrary analytic curves and normal frames.

Let $\tilde c(t)$ be an analytic curve with analytic normal frame
 $\tilde c'(t), \tilde n_1(t), \tilde n_2(t)$.
Assume that
 that $|c(t)-\tilde c(t)|\le Ct^2$ and $|n_j(t)-\tilde n_j(t)|\le C t^2$ for
 $t\in(-\epsilon,\epsilon) $ and $j=1,2$. Here we think of $c(t)$ as being an osculating circle for $\tilde c(t)$ at $\tilde c(0)$, but our argument below works for any curves that are close to second order.

We assume that both $c$ and $\tilde c$ are parameterized by arc length.
 Introduce the  spinning normal fields for $\tilde c$
\[
\tilde n(t) = \cos(a t)\tilde n_1(t) + \sin(at)\tilde  n_2(t).
\]

Define the Bj\"orling surfaces
\[
 \tilde F (z) = \re\left( \tilde c(z) -  i \int_0^z \tilde n(w)\times \tilde c'(w)\, dw \right).
\]

As $c'\times n_1=n_2$, we have
\[
c'\times n =  \cos(a t) n_2(t)- \sin(at)  n_1(t)
\]
(and similarly for $\tilde c$).

This allows us to estimate the distance between the parameterizations $F(z)$ and $\tilde F(z)$.

\begin{lemma}
\label{lem:comp}
For $|\im z|\le d =\frac{\log (a)}{a}$ and $|\re z|<\frac{\pi}{a}$ and sufficiently large $a$
\[
|\tilde F(z)-F(z)|\le 6C \frac{(\log a)^2}{a^2}
\]
\end{lemma}
\begin{proof}
The idea is to integrate
 by parts twice.
\begin{align*}
|\tilde F(z)-F(z)|
&\le \left|
\tilde c(z)-c(z) - i\int_0^z  \cos(a w) (\tilde n_2(w)-n_2(w))- \sin(a w) 
 (\tilde n_1(w)-n_1(w)) \, dw
\right|
\\
&=C|z|^2
+\left|
\frac{\sin(a w)}{a}(\tilde n_2(w)-n_2(w))- \frac{\cos(a w)}{a}
 (\tilde n_1(w)-n_1(w))
 +\right|_{w=0}^z +
\\&\qquad+\left|\int_0^z  \frac{\sin(a w)}{a} (\tilde n_2'(w)-n_2'(w))-
 \frac{\cos(a w)}{a}  (\tilde n_1'(w)-n_1'(w)) \, dw
\right|
\\
&\le
C|z|^2(1+|\frac{\sin(a w)}{a}|+|\frac{\cos(a w)}{a}|)
+C|z|(|\frac{\sin(a w)}{a^2}|+|\frac{\cos(a w)}{a^2}|)+
\\&\qquad
+\left|
\int_0^z  \frac{\cos(a w)}{a^2} (\tilde n_2''(w)-n_2''(w))-
 \frac{\sin(a w)}{a^2}  (\tilde n_1''(w)-n_1''(w)) \, dw
\right|
\\
&\le
C|z|^2(1+|\frac{\sin(a w)}{a}|+|\frac{\cos(a w)}{a}|)
+2C|z|(|\frac{\sin(a w)}{a^2}|+|\frac{\cos(a w)}{a^2}|).
\end{align*}

Now we use that  the domain is $|\im z|\le d =\frac{\log (a)}{a}\approx \frac{\sinh^{-1}(a)}{a}$
 and $|\re z|<\frac{\pi}{a}<\frac{\log (a)}{a}$, so that for $a$ large, $|z|<\sqrt{2}d$. In this
 domain, we get
\[
|\tilde F(z)-F(z)|\le 6Cd^2+2C \frac{d}{a}.
\]
As $d= \frac{\log (a)}{a}$,
\[
|\tilde F(z)-F(z)|\le 6C \frac{(\log a)^2}{a^2}.
\]
\end{proof}

\vspace{.2cm}

\section{The proof of Theorem~\ref{thm1} in the circular case.}

In this section, we  prove the following  version of
 Theorem~\ref{thm1} in the circular case.

Let $T_R$ be the
 $(R - \frac{1}{R})$-neighborhood in $\rth$ of the circle
 $\esf^1(R) = \{ x_1^2 + x_2^2 = R^2 \} \subset
\R^2 \times \{ 0 \}$. This domain is the region within which we want to consider embeddedness first.

\begin{theorem} \label{thm2}
For each $a \in \N$ and $R > 1$, let $\widehat{H}_{a,R}$ be the
 component of the embedded minimal disk $F([\frac{- \pi}{2a},
 \frac{\pi}{2a}] \times
\R) \cap T_R$ containing the circular arc $F([\frac{- \pi}{2a},
\frac{ \pi}{2a}] \times \{ 0 \})$. For $a \in \N$ sufficiently large, the orbit $H_{a,R}$ of
 $\widehat{H}_{a,R}$ under the cyclic group $\Z_{2a}$ generated by 
rotation around the positive $x_3$-axis by the angle $\frac{\pi}{a}$,
 is an embedded minimal annulus. Furthermore, as $a \to \infty$  and $R$ is fixed, the
 annuli $H_{a,R}$ converge to the minimal foliation $\lc_R$ of $T_R$ 
consisting of the flat
 disks of radii $R- \frac {1}{R}$ centered at points of $\esf^1(R)$ and orthogonal
 to $\esf^1(R)$ and with $S(\lc_R) = \esf^1(R)$. In particular, there exists a
 divergent sequence of $R_n \to \infty$ as $n \to \infty$, so that the 
bent helicoids $\overline{H}_n = H_{n,R_{n}}$ are embedded and converge to the
 Colding-Minicozzi limit minimal lamination $\lc$ of $\rth - x_3$-axis 
consisting of leaves which are half planes with axis the $x_3$-axis and with $S(\lc) = \esf^1(1)$.
\end{theorem}

\begin{figure}[H]
  \centerline{\includegraphics[width=3in]{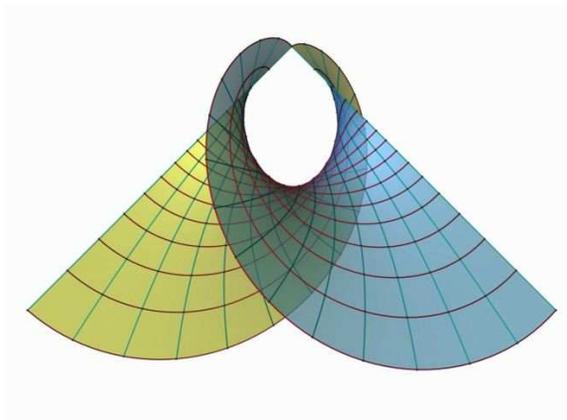}}
  \caption{A fundamental piece of $H_a$ in a ``sector'' for $a=2$.}
  \label{fig:embedded2}
\end{figure}

\begin{proof}
Fix $R > 1$ and $a \in \N$. 
In Corollary 1 in section 3, we proved that the
 $\widehat{H}_a = F((\frac{-\pi}{2a}, \frac{\pi}{2a}]
 \times (-\infty, \infty)$ is a fundamental piece of 
the surface $H_a$, which is embedded in $\rth$ with boundary
 two straight lines. We now denote these two lines by $L(\frac{-\pi}{2a}),
 L(\frac{\pi}{2a})$ and remark that they lie in the
 $(x_1, x_2)$-plane and make
 an angle of $\frac{\pi}{a}$ at the origin. Let $H_{a,R}$
 denote the immersed minimal surface component of $H_a \cap T_R$
 that contains $\esf^1(1)$. Note that $H_{a,R}$ is the $\Z_{2a}$-orbit of the embedded component $\widehat{H}_{a,R}$ of
 $\widehat{H}_a \cap T_R$ that contains the circular arc
 $\sigma_a = F([\frac{-\pi}{2a}, \frac{\pi}{2a}] \times \{ 0 \})$,
 where $\Z_{2a}$ is generated by rotation by $\frac{\pi}{a}$ around the
 $x_3$-axis. More precisely, $H_{a,R}$ is the image under $F$ of the component
 of $F^{-1}(T_{R}) \cap ([ -\frac{\pi}{2a}, \frac{\pi}{2a}] \times
(-\infty, \infty))$ that contains the interval $[ -\frac{\pi}{2a}, \frac{\pi}{2a}] \times \{ 0 \}$. 

By the results in the previous section, for large values of $a \in \N$, 
the surface $\widehat{H}_{a,R}$ is a compact embedded disk which intersects
 $\partial T_R$ almost orthogonally in two almost circular arcs in $\partial T_{R}$.  These arcs join
 the end points of line segments $l(\frac{-\pi}{2a})\subset
L(\frac{-\pi}{2a}) \cap
T_R$, $ l(\frac{\pi}{2a}) \subset L(\frac{\pi}{2a})
\cap T_R$, which make up the remainder of
 $\partial \widehat{H}_{a,R}$. If the embedded disk 
$\widehat{H}_{a,R}$ were contained in the sector of $\rth$
 containing the circular arc $\sigma_a$ and bounded by the vertical half planes containing
 $l(\frac{-\pi}{a}), l(\frac{\pi}{2a})$, respectively, then the
 $\Z_{2a}$-orbit $H_{a,R}$ of $\widehat{H}_{a,R}$ would be an embedded
 annulus. Although  $\widehat{H}_{a,R}$ fails to be contained in this
 sector (see Figure \ref{fig:embedded2}), we shall still be able to prove that $H_{a,R}$ is an embedded minimal annulus
 for $a\in \N$ large.
 
Consider $T_R$ with the ``cylindrical'' coordinates $(\theta, x)$ with
 $x \in D_R$, where $D_R$ is the disk in the $(x_1, x_3)$-plane of
 radius $R- \frac{1}{R}$ centered at the point $(R, 0, 0)$. For $a \in \N$ sufficiently large and for $\ve(a) = \frac{1}{a}$, the sequence of $\ve(a)$-tubular neighborhoods $\esf^1_{\ve(a)}(1)$ of $\esf^1(1)$ in $T_R$, when intersected with $H_{a,R}$ and then translated by $(-1, 0, 0)$ and expanded homothetically by 
the factor $a$, produces a sequence of minimal annuli which converges on compact subsets of $\rth$ to a helicoid intersected with the solid cylinder of radius $1$ around the $x_2$-axis (this follows by direct calculation).
For $x \in D_{R}$, let $C_x$ denote the horizontal 
circle $\{ (\theta, x) \mid 0
 \leq \theta < 2\pi \}$ in our cylindrical coordinates of $T_R$ in $\rth - x_3$-axis. It follows that for $a\in \N$ large, $\widehat{H}_{a,R}$ intersects every horizontal circle
 $C_x \subset \esf^1_{\ve(a)} (1)- \esf^{1}(1)$,
 transversely in a single point.
 Furthermore, for $a$ large, every horizontal circle
 $C_x \subset (T_R - \esf^1_{\ve(a)}(1))$, intersects $\widehat{H}_{a,R}$ transversely 
in a single point and the angle of intersection is uniformly bounded
 away from zero by a positive constant which is independent
 of $a$ (this follows
 from our formula for the Gauss map of $\widehat{H}_{a,R}$ and the
 estimates in the previous section). Since every such circle $C_x \subset (T_{R} - \esf^1(1))$ is invariant
 under $\Z_{2a}$, it follows that for $a \in \N$ sufficiently large, 
$H_{a,R}$ is an embedded minimal annulus.

 In cylindrical
 coordinates, we see that for $a \in \N$ sufficiently large,
 $H_{a,R} - \esf^1(1)$ is a two component multigraph over $D_R - \{(1,0,0) \}$
 invariant under the action of $\Z_{2a}$. Also note that each of these
 multigraphs is stable with a positive Jacobi function induced by
 the killing field of $\rth$ generated by rotation around the $x_3$-axis.
 In particular, by curvature estimates for stable minimal surfaces \cite{sc3},
 we see that the sequence of surfaces $\{ H_{a,R} \}_{a \in \N}$ has
 uniformly locally bounded curvature in any ball in $T_R$ of positive 
distance from $\esf^1(1)$. 

It is now standard (e.g. see the proof of
 Theorem 1.6 in \cite{mr8}) that a subsequence of the surfaces 
$\{ H_{a,R} - \esf^{1}(1) \}_{a \in \N}$ converges $C^2$ to a minimal lamination $\widehat{\lc}_R$ of
 $T_R - \esf^1(1)$ whose leaves are mapped to other leaves under any
 rotation around the $x_3$-axis. Clearly, the leaves of $\widehat{\lc}_R$ are
 punctured flat disks of radius $R- \frac{1}{R}$ which are orthogonal to
 $\esf^1(1)$ and are centered along $\esf^1(R)$ (for example, consider
 the values of the Gauss map of $H_{a,R}$ away from $\esf^1(1)$ for $a\in \N$ large). Thus, $\widehat{\lc}_R$ 
extends to the foliation $\lc_R$ of $T_R$ by flat disks orthogonal to
 $\esf^1(1)$. Since the tangent planes of the $H_{a,R}$ are orthogonal
 to the tangent planes to $\lc_R$ along  $ \esf^1(1)$, 
the sequence $\{ H_{a,R} \}_{a\in \N}$ converges to $\lc_R$ with 
singular set of $C^{1}$-convergence $S(\lc_R) = \esf^1(1)$. This 
concludes the proof of the first statement of the theorem. The second
 statement then follows by applying a standard diagonal argument.
\end{proof}

\section{The proof of Theorem~\ref{thm1} in the analytic case.}

In the last section, we proved Theorem~\ref{thm1} in the case the curve
 $\alpha$ is the unit circle in the $(x_1, x_2)$-plane and our open set is $\rth$. We now prove the
 theorem in the special case where $\alpha$ is a properly embedded
 analytic curve in an open set $O$.

In what follows, it suffices  to consider $\alpha$ an open curve. 
If the curve is closed, one faces the additional technical complication that the normal fields need to close up.
 Without loss of generality, we may assume that our analytic curve $\alpha(t)$ has unit speed with anlaytic frame $\alpha'(t), n_1(t), n_2(t) =
\alpha'(t) \times n_1(t)$. Let 
$$ n(t) = \cos(at) n_1(t) +
\sin (at) n_2(t).$$
For $a \in \R^+$, let $H_{\alpha,a}$ be related Bj\"orling surface or {\em bent helicoid}. Fix a point $p \in \alpha(t)$, which we may assume has the form $p= \alpha(0)$. After a rigid motion, we may assume that 

$$\alpha'(0) = (1,0,0), \; \langle \alpha''(0),
(-1,0,0) \rangle = \kappa \geq 0.$$

\noindent If the curvature $\kappa \neq 0$, then, after a dilation, we also 
may assume that $\kappa =1$. 

The analysis of the case $\kappa = 0$ 
and the case $\kappa=1$ are similar. We  only consider the
 case $\kappa \neq 0$; in both cases, one compares the geometry of 
$H_{\alpha,a}$ with the geometry of a standard surface where the standard 
surface is a helicoid if $\kappa =0$ or the bent helicoid $H_a$ when $\kappa = 1$. So, assume now that $\kappa =1$. 

Consider a continuous family 
${\cal D}$ of pairwise disjoint disks $D_t$ which are orthogonal
 to $\alpha(t)$ for each $t$ and which lie in the interior of another 
such family $\widetilde{{\cal D}}$. Note that for $R > 1$ fixed
 and large, there exists a small $\ve > 0$ such that the family of disks $D_t, - \frac{\ve}{2} \leq
t \leq \frac{\ve}{2}$, are embedded and contained in the domain $T_R(\ve) = \{ (\theta, x) \in T_R \mid -\ve \leq
\theta \leq \ve\}$, where $T_R$ is defined just before the statement of Theorem~\ref{thm2} and the cylindrical coordinates on $T_R$ are those introduced in the proof of Theorem~\ref{thm2}. Theorem~\ref{thm1} in the case $\alpha$ is analytic easily follows from the following assertion, after restricting neighborhoods appropriately.

\begin{assertion}
Fix $R >1$. Then there exists a small $\ve > 0$ such that for $a \in \N$ sufficiently large, the component $H_{\alpha,a}(\ve)$ of 
$F_a( [-\ve, \ve] \times (-\delta(\ve),
\delta(\ve)) \cap T_R (2\ve)$ containing $\alpha_{\ve} = F_a([-\ve, \ve] \times \{ 0\})$ is an embedded disk. Here, the domain $[-\ve, \ve] \times (-\delta(\ve), \delta(\ve))$ is a box neighborhood of $[-\ve, \ve] \times \{0 \} \subset \C$, where the Bj\"orling data is defined. 
\end{assertion}

\begin{proof}
We  first consider the special case where $n_1(0) =(-1, 0, 0)$. Let $\widetilde{\alpha}(t) = (\cos(t), \sin(t), 0)$ and
$\widetilde{n}_1(t) = -\widetilde{\alpha}(t)$, 
$\widetilde{n}_2(t)= (0, 0, 1)$. The data for $\alpha$ and $\widetilde{\alpha}$ agree to second order at $t=0$ and the data for $\widetilde{\alpha}$  produces the Bj\"orling bent helicoids $H_n$ for $n \in \N$. By the proof of Theorem~\ref{thm2}, for $a \in \N$ large, the component $H_{\widetilde{\alpha}, a}( \ve)$ is an embedded disk. 

As in the case $\alpha$ was the circle $\widetilde{\alpha}$, which we considered in the previous section, for $a \in \N$ large, on the scale of curvature 
and 
 for any sufficiently small $\ve >0$, the surface $H_{\alpha, a}(\ve)$ is 
closely approximated near each point of $\alpha$ by homothetically shrunk helicoids initially contained in a cylinder of radius $1$ around its axis, in the following sense. At every point $q \in \alpha_{\ve} = \alpha ([-\ve, \ve])$ and inside 
the $\frac{1}{a}$-neighborhood 
$N_{\rth}(\alpha_{\ve}, \frac{1}{a})$ of $\alpha_{\ve}$ in $\rth$ 
the related sequence of surfaces under dilations by the factor 
$a$ at $q$, converge to a helicoid with axis tangent to $\alpha$ 
at $q$. Note that
 $N_{\rth}(\alpha_{\ve}, \frac{1}{a}) \cap H_{\alpha, a}( \ve) =
N(\alpha_{\ve}, \frac{1}{a})$ is a simply connected neighborhood of 
$\alpha_{\ve}$ in $H_{\alpha, a}(\ve)$. Using the fact that $H_{\alpha, a}(\ve)$ is closely approximated by
 a ruled surface, shows that for $a \in \N$ sufficiently large, the self-intersection 
set of $H_{\alpha, a}( \ve)$ is disjoint from the 
$\frac{1}{a}$-neighborhood of $\alpha_{\ve}$.

Let $\eta=\frac{2\pi}{a}$.
Our
 approximation results imply that for $a \in \N$ large that 
$H_{\alpha, a}( \eta) - N_{\rth}
(\alpha_{\eta}, \frac{1}{a})$ consist of two parametrized disks 
$S_+, S_{-}$ that are multi-graphs over $H_{\widetilde{\alpha}, a}( \eta) -
N_{\rth}(\alpha_{\eta}, \frac{1}{a})$ of norm on the order of 
$(\frac{\log(a)}{a})^2$ for $a$ sufficiently large (this estimate also depends on $R$ but since $R$ is fixed it can be assumed to be uniform in $a$).
 On the other hand, for $a$ large, the distance between 
successive sheets of the two spiraling multigraphs $\widetilde{S}_{+},
 \widetilde{S}_{-} \subset (H_{\widetilde{\alpha}, a}( \eta)- N_{\rth}
(\alpha_{\eta}, 
\frac{1}{a}))$ is bounded from below by $\frac{C'}{a}$, where $C'$ 
is a positive constant depending only on $R$.

When choosing $a \in \N$ large, 
the sheets of $\widetilde{S}_{+}$ separate the sheets of the multigraphs
 $S_{-}$ from each other (similarly $\widetilde{S}_{-}$ separates the 
sheets $S_{+}$ from each other), then the sheets of $S_{+}$ and $S_{-}$ do not intersect. Since $S_+$ can be expressed 
as a small graph over $\widetilde{S}_+$ with gradient bounded 
uniformly for $a \in \N$ large. Hence, $H_{\alpha, a}(\ve)$
 is an embedded disk for some fixed small $\ve > 0$. 

This completes the proof of the assertion under the assumption that $n_1(0) = (-1, 0, 0)$. In the case, $n_1(0) \neq (-1, 0, 0)$, one compares the surface
$H_{\alpha, a}( \ve)$, for large $a \in \N$, with $\overline{H}_a(\theta)$ where
$\overline{H}_a(\theta)$ is the bent helicoid $\overline{H}_{a}$ rotated 
so that the normal fields satisfy:
$$n_1(a, \theta)(0) = n_1(0).$$
Then one proceeds as above. Thus, the general case follows from our special case where $n_1(0)= (-1, 0, 0)$. This completes the proof of the assertion.

\end{proof}

\section{The proof of Theorem~\ref{thm1} in the $C^{1,1}$-case.}

Consider now an arbitrary properly embedded $C^{1,1}$-curve $\alpha$
 in an open set $O$ of $\rth$. We  can just consider the case 
where $\alpha$ is noncompact because the compact case follows from the
 same arguments. Consider $\alpha \colon (a,b) \to \rth$ to be a unit
 speed $C^{1,1}$-parametrization of the image curve $\alpha$. Fix a
 compact exhaustion
$$[a_1, b_1] \subset ... \subset
[a_n, b_n] \subset ...$$
of $(a,b)$. Recall that a $C^{1,1}$-curve $\alpha(t)$ has locally bounded
 curvature function $\kappa(t)$ defined almost everywhere; in fact, $\alpha'(t)$ is absolutely continuous with a related locally bounded difference 
quotient function $\widehat{\kappa}(t)$. Since $\alpha$
 is a $C^{1,1}$-curve, 
there exists a sequence of embedded unit speed analytic curves
 $\beta_n \colon 
[a_n, b_n ] \to \rth$ which converge $C^1$ to $\alpha$. The $\beta_n$ can be chosen to  have
 uniformly bounded curvature at most $\min \widehat{\kappa}|_{[a_{k},
b_{k}]}$ on any fixed $[a_k, b_k] 
\subset [ a_n, b_n]$ for $n \geq k$. Their related curvature
 functions are uniformly bounded by $\widehat{\kappa}(t)$.
 
 To see this, first
 convolve the  $C^{1,1}$-curve with a mollifier. This gives a $C^\infty$-curve which will be uniformly close to the original curve. The Lipschitz bound on the velocity then bounds the second derivative of the mollified curve. (Differentiate the convolution twice and integrate by parts once).
 Then, these $C^\infty$-curves can be approximated by analytic curves $\beta_n$ converging to $\alpha$ and with curvature functions converging to $\widehat{\kappa}(t)$.
 It then follows 
from arguments of the previous
 section that for fixed $n$ and $k$ with $n > k$, there is a sequence
 of bent helicoids containing $\beta_n [a_k, b_k]$ which give rise to
 a Colding-Minicozzi minimal lamination of the $\lambda \widehat{\kappa}(t)$-normal bundle of $\beta_n [a_k, b_k]$ for any positive $\lambda < 1$. A standard diagonal argument together with arguments from the analytic case
 then produces a sequence of bent helicoids that converges to a limit 
 minimal lamination satisfying the requirements of Theorem~\ref{thm1}.

\center{William H. Meeks, III at bill@gang.umass.edu,\\
Math Department, University of Massachusetts, Amherst, MA  01003}.\\

\center{Matthias Weber matweber@indiana.edu,\\
Math Department, University of Indiana, Bloomingtion, IN  47405}.

\bibliographystyle{plain}

\bibliography{bill1}

\begin{thebibliography}{10}

\bibitem{cm23}
T.~H. Colding and W.~P. Minicozzi~II.
\newblock The space of embedded minimal surfaces of fixed genus in a 3-manifold
  {I}{V}; {L}ocally simply-connected.
\newblock {\em Annals of Math.}, 160:573--615, 2004.

\bibitem{dhkw}
U.~Dierkes, S.~Hildebrandt, A.~K\"{u}ster, and O.~Wohlrab.
\newblock {\em Minimal Surfaces I}.
\newblock Grundlehren der mathematischen {W}issenschaften 296. Springer-Verlag,
  1992.

\bibitem{eww1}
T.~Ekholm, B.~White, and D.~Wienholtz.
\newblock Embeddedness of minimal surfaces with total curvature at most $4\pi$.
\newblock {\em Annals of Math.}, 155:209--234, 2002.

\bibitem{me7}
W.~H. Meeks~III.
\newblock The classification of complete minimal surfaces with total curvature
  greater than {$-8\pi$}.
\newblock {\em Duke Math. J.}, 48:523--535, 1981.

\bibitem{me25}
W.~H. Meeks~III.
\newblock The regularity of the singular set in the {C}olding and {M}inicozzi
  lamination theorem.
\newblock {\em Duke Math. J.}, 123(2):329--334, 2004.

\bibitem{me30}
W.~H. Meeks~III.
\newblock The lamination metric for a {C}olding-{M}inicozzi minimal lamination.
\newblock {\em Illinois J. of Math.}, 49:645--658, 2005.

\bibitem{mpr10}
W.~H. Meeks~III, J.~P\'{e}rez, and A.~Ros.
\newblock Embedded minimal surfaces: removable singularities, local pictures
  and parking garage structures, the dynamics of dilation invariant collections
  and the characterization of examples of quadratic curvature decay.
\newblock Preprint.

\bibitem{mpr8}
W.~H. Meeks~III, J.~P\'{e}rez, and A.~Ros.
\newblock The geometry of minimal surfaces of finite genus {I}{I}{I}; bounds on
  the topology and index of classical minimal surfaces.
\newblock Preprint.

\bibitem{mpr3}
W.~H. Meeks~III, J.~P\'{e}rez, and A.~Ros.
\newblock The geometry of minimal surfaces of finite genus {I}; curvature
  estimates and quasiperiodicity.
\newblock {\em J. of Differential Geometry}, 66:1--45, 2004.

\bibitem{mpr4}
W.~H. Meeks~III, J.~P\'{e}rez, and A.~Ros.
\newblock The geometry of minimal surfaces of finite genus {I}{I}; nonexistence
  of one limit end examples.
\newblock {\em Invent. Math.}, 158:323--341, 2004.

\bibitem{mr13}
W.~H. Meeks~III and H.~Rosenberg.
\newblock The minimal lamination closure theorem.
\newblock {\em Duke Math. J. {\em (to appear)}}.

\bibitem{mr10}
W.~H. Meeks~III and H.~Rosenberg.
\newblock The theory of minimal surfaces in ${M} \times \mathbb{R}$.
\newblock {\em Comment. Math. Helv. {\em (to appear)}}.

\bibitem{mr8}
W.~H. Meeks~III and H.~Rosenberg.
\newblock The uniqueness of the helicoid and the asymptotic geometry of
  properly embedded minimal surfaces with finite topology.
\newblock {\em Annals of Math.}, 161:727--758, 2005.

\bibitem{sc3}
R.~Schoen.
\newblock {\em Estimates for Stable Minimal Surfaces in Three Dimensional
  Manifolds}, volume 103 of {\em Annals of Math. Studies}.
\newblock Princeton University Press, 1983.
\newblock MR0795231, Zbl 532.53042.

\end{thebibliography}

\end{document}